\documentclass[a4paper,10pt]{article}
\usepackage[american]{babel}
\usepackage{amssymb,amsmath,amsthm}
\newtheorem{theorem}{Theorem}
\newtheorem{corollary}[theorem]{Corollary}
\newtheorem{example}{Example}
\newcommand{\Nat}{\mathbb{N}}
\newcommand{\Real}{\mathbb{R}}
\newcommand{\set}[1]{\left\{#1\right\}}
\newcommand{\ds}[1]{_{_{#1}}}

\begin{document}
\title{An exact expression of $\pi(x)$}

\author{ J. E. Palomar Taranc\'{o}n}

\maketitle

\begin{abstract}
The author states an exact expression denoting the distribution of primes.
\end{abstract}
\section{Introduction}
Although the distribution of primes $\pi(x)$ is only approached by expressions like $x/ln(x)$, we state an exact expression
of $\pi(x)$. Such an expression of $\pi(x)$ is piecewise defined in intervals of the form $(p\ds{n},p\ds{n+1}^2)$, hence by recursion over the previous intervals the whole function $\pi(x)$  can be built.
The expression arises from the Theorem~\ref{thm1} which can be obtained easily using the number representation defined in
\cite{palomar} as an instance of the abstraction concept introduced by the author.

\section{Distribution of primes}
Let $\mathbf{P}=\set{p\ds{1},p\ds{2}\dots{}p\ds{n}\dots}=\set{2,3,5,7\dots}$ stand for the ordered sequence of all primes, and for every
couple of positive integers $n,m\in\Nat$  such that $m\leq{}n$, let $\sigma\ds{n,m}(x)$, denote the  function defined by
\begin{equation}
 \label{eq:sym}
\forall x\in\Real:\qquad \sigma\ds{n,m}(x) =\sum_{1\leq{}k_1<k_2\dots<k_m\leq{}n}\left[\frac{x}{p\ds{k_1}p\ds{k_2}\dots{}p\ds{k_m}}\right]
\end{equation}
where $[\,\,]:\Real\rightarrow\mathbb{Z}$ stands for the floor function.

From these  functions define $\gamma\ds{n,m}(x)$ recursively, as follows.
\begin{align*}
 &\gamma\ds{n,n}(x)= \sigma\ds{n,n}(x)\\
&\gamma\ds{n,(n-1)}(x)=\sigma\ds{n,(n-1)}(x)-\binom{n}{n-1}\gamma\ds{n,n}(x)\\
&\gamma\ds{n,(n-2)}(x)=\sigma\ds{n,(n-2)}(x)-\binom{n-1}{n-2}\gamma\ds{n,(n-1)}(x)-\binom{n}{n-2}\gamma\ds{n,n}(x)\\
&\dots\dots\dots\dots\dots\dots\dots\dots\dots\dots\dots\dots\dots\dots\dots\dots\dots\dots\\
&\gamma\ds{n,m}(x)=\sigma\ds{n,m}(x)-\sum_{k=m+1}^n\binom{k}{m}\gamma\ds{n,k}(x)\\
&\dots\dots\dots\dots\dots\dots\dots\dots\dots\dots\dots\dots\dots\dots\dots\dots\dots\dots
\end{align*}
Finally, let $\Upsilon_n(x)$ be the function,
\begin{equation}
 \label{eq:psi}
\Upsilon_n(x)=[x]-\sigma\ds{n,1}(x)+\sum_{k=2}^n(k-1)\cdot\gamma\ds{n,k}(x)+n-1
\end{equation}
\begin{theorem}\label{thm1}
 For each $n>1$ and for every $x\in(p\ds{n},p^2\ds{n+1})$ the following equality holds.
\begin{equation}
 \label{eq:thm}
\Upsilon_n(x)=\pi(x)
\end{equation}
\end{theorem}
\begin{corollary}
 For every $x\geq0$,
\begin{equation}
 \pi(x)=
\begin{cases}
 0 \text{ if \/} x<2\\
1 \text{ if  \/} 2\leq x<3\\
2 \text{ if \/} x=3\\
\Upsilon_n(x)\text{ if } p\ds{n}<x<p^2\ds{n+1}\text{ and } n\geq 2
\end{cases}
\end{equation}
\end{corollary}

Indeed, because $\pi(x)$ detects primes, by the former definition the function $\pi(x)$ can be built recursively.
\begin{example}
 To calculate $\pi(x)$ in the interval $(p\ds{4},p\ds{5}^2)=(7,121)$ we use the function $\Upsilon_4(x)$, which can
be written as follows.
\begin{multline}
 \Upsilon_4(x)=[x] -\left(\left[\frac{x}{2}\right]+\left[\frac{x}{3}\right]+\left[\frac{x}{5}\right]+\left[\frac{x}{7}\right]\right)+
\left(\left[\frac{x}{2\cdot3}\right]+\left[\frac{x}{2\cdot5}\right]+\left[\frac{x}{2\cdot7}\right]+\right.\\
\left.\left[\frac{x}{3\cdot5}\right]+\left[\frac{x}{3\cdot7}\right]+\left[\frac{x}{5\cdot7}\right]\right)-\left(\left[\frac{x}{2\cdot3\cdot5}\right]+
\left[\frac{x}{2\cdot3\cdot7}\right]+\left[\frac{x}{2\cdot5\cdot7}\right]+\left[\frac{x}{3\cdot5\cdot7}\right]\right)-\\
11\left[\frac{x}{2\cdot3\cdot5\cdot7}\right]+3
\end{multline}
and it is not difficult to see, that for every $x\in(7,\,121)$:\quad $\Upsilon_4(x)=\pi(x)$.
Of course, for the interval $(11,169)$ one can use the function $\Upsilon_5(x)$, for $(13,289)$ the function
$\Upsilon_6(x)$ and so on.
\end{example}

\end{document}